\input amstex\documentstyle{amsppt}  
\pagewidth{12.5cm}\pageheight{19cm}\magnification\magstep1
\topmatter
\title On the group attached to a special
Weyl group representation
\endtitle
\author G. Lusztig\endauthor
\address{Department of Mathematics, M.I.T., Cambridge, MA 02139}\endaddress
\thanks{Supported by NSF grant DMS-2153741}\endthanks
\endtopmatter   
\document

\define\Irr{\text{\rm Irr}}
\define\sneq{\subsetneqq}

\define\op{\oplus}
   
\define\part{\partial}

\define\m{\mapsto}
\define\do{\dots}

\define\sub{\subset}    
\define\bxt{\boxtimes}
\define\T{\times}

\define\nl{\newline}
\redefine\i{^{-1}}

\define\ot{\otimes}

\define\sg{\text{\rm sgn}}

\redefine\d{\delta}
\define\e{\epsilon}

\redefine\G{\Gamma}
\redefine\D{\Delta}

\define\CC{\bold C}

\define\NN{\bold N}

\define\ZZ{\bold Z}

\define\sha{\sharp}

\subhead 1.1\endsubhead
Let $W$ be a Weyl group. Let $\Irr(W)$ be the set of (isomorphism
classes of) irreducible representations over $\CC$ of $W$.
Let $\Irr_{sp}(W)$ be the subset of $\Irr(W)$ consisting of special representations.
In \cite{L84, 13.1} to each $E\in\Irr_{sp}(W)$ we have
associated a certain finite group; we denote this finite
group by $\G_E$.
The definition of $\G_E$ in {\it loc.cit.} was in algebraic
geometric terms using Springer representations.
In this paper we will associate to $E$ a finite group
$\D_E$ (isomorphic to $\G_E$) using results of \cite{L25}; thus our definition is based on pure algebra (not geometry). 

\subhead 1.2\endsubhead
When $W$ is a product $W_1\T W_2\T\do\T W_e$ where $W_i$ are
irreducible Weyl groups and $E=E_1\bxt E_2\bxt\do\bxt E_e$ where
$E_i\in \Irr_{sp}(W_i)$ for $i=1,2,\do,e$, then
$\D_E=\D_{E_1}\T\D_{E_2}\T\do\T\D_{E_e}$. Thus it is enough to
define $\D_E$ when $W$ is irreducible.
If $E\in\Irr_{sp}(W)$ is not superspecial in the sense of
\cite{L25}, then for some
$I'\sneq I$ and some $E'\in\Irr_{sp}(W_{I'}$
we have that $J_{W_{I'}}^W(E')=E$
or $J_{W_{I'}}^W(E')\ot\sg_W=E$
(here $\sg_W$ is the sign representation of $W$);
moreover we have $\D_E\cong\D_{E'}$. Thus it is enough to
define $\D_E$ when $W$ is irreducible and $E$ is
superspecial.

For any $E\in\Irr(W)$ let $c_E\in\NN$ be as in
\cite{L25,0.1}.
We shall prove

\proclaim{Proposition 1.3} Assume that $W$ is irreducible
and $E\in\Irr(W)$ is superspecial. There is a
unique finite group $\D_E$ (up to isomorphism) which has
properties (a)-(d) below.

(a) $\sha(\D_E)=c_E$.

(b) For any $g\in\D_E$, $g,g\i$ are conjugate.

(c) The list of sizes of the various conjugacy classes of
$\D_E$ is the same as the list
$\{c_E/c_{E_j};j=1,2,\do,t\}$ where
$E_1\op E_2\op\do\op E_t$ is the decomposition of $Z_W$
(in \cite{L25,\S2}) in irreducible representations of $W$.
(See 2.2(a).)

(d) If $c_E>1$, then $\D_E$ has a subgroup $\D'_E$ of
index $2$.
\endproclaim
(d) can be refined as follows.
The list of sizes of the various conjugacy
classes of $\D_E$ decomposes into a union of two lists,
one involving conjugacy classes contained in $\D'_E$, the
other one involving conjugacy classes contained in
$\D_E-\D'_E$. Under (c) this corresponds to a partition
of $\{1,2,...,t\}$ in which one part is
$\{j;b_{E_j}\text{even}\}$ and one part is
$\{j;b_{E_j}\text{odd}\}$. 
(Here $b_{E_j}$ are as in \cite{L84,(4.1.2)}.
See 2.2(a1).) This gives actually a canonical choice for
$\D'_E$. (This refinement of (d) will not be used.)

\subhead 1.4\endsubhead
Note that $c_E$ is one of the numbers

$1,2,2^2,2^3,\do$, $3!,4!,5!$.
\nl
If $\D_E$ is

$\{1\},S_2,S_2\T S_2,S_2\T S_2\T S_2,...$, $S_3,S_4,S_5$
\nl
(respectively) then properties
1.3(a)-(d)
are satisfied. (Here $S_n$ is the symmetric group in
$n$ letters.) In the rest of the proof
we assume that $\D=\D_E$ is a group satisfying 1.3(a)-(d).
We will show that $\D$ is uniquely dermined up to
isomorphism.
If $c_E$ is $1$ or $2$ then $\D$ must be $\{1\}$ or $S_2$
by 1.3(a).
If $c_E=2^k$ for some $k\ge2$ then by 1.3(a) we have
$\sha(\D)=2^k$ and by 1.3(c), all
conjugacy classes in $\D$ have size equal to $2^k/2^k=1$
so that $\D$ is abelian. Using 1.3(b) we see that
any $g\in\D$ satisfies $g^2=1$. It follows that
$\D$ is a product of $k$ copies of $S_2$ as desired.
(Note that
1.3(d) has not been used.)

In the remainder of the proof we assume that $c_E$ is one of
$3!,4!,5!$. We will associate to $\D$ a set $S$ consisting
of $3,4,5$ elements. An element of $S$ will be

one of the three elements of order $2$ in $\D$ (if
$c_E=3!$);

one of the four pairs of elements of order $3$ one being
the inverse of the other (if $c_E=4!$);

one of the five triples of commuting elements of order $2$
each of which has centralizer of order $8$ (if $c_E=5!$).

Let $Perm(S)$ be the group of permutations of $S$. Then
$\D$ acts naturally on $S$; this gives a homomorphism
$\D@>>>Perm(S)$. We will show that it is an isomorphism.

Assume first that $c_E=3!$. The conjugacy classes
in $\D$ are $A_6^1,A_3^2,A_2^3$.
(We shall usually write $A_\e^\d$ for a conjugacy class in a
finite group which consists of $\d$ elements and $\e$ is
the order of the centalizer $Z(z)$ of an element $z$ in the
conjugacy class.) Let $S=A_2^3$. We have $\sha S=3$.
Let $s\in S$. We have $\sha Z(s)=2$, 
hence $Z(s)=\{1,s\}$. Let $K$ be the set of all $g\in\D$
such that $gsg\i=s$ for all $s\in S$. Then $g\in\{1,s\}$
for any $s\in S$ hence $g=1$. Thus the homomorphism
$\D@>>>Perm(S)$, $g\m(s\m gsg\i)$ with kernel $K$ is
injective. Since $\sha(\D)=\sha(Perm(S))$ we see that the
last homomorphism is bijective so that $\D\cong S_3$ as desired. (We have not used 1.3(b),(d).)

We now assume that $c_E=4!$. The conjugacy classes
in $\D$ are $A_{24}^1,A_8^3,A_4^6,{}'A_4^6,A_3^8$.

Let $u\in A_3^8$. Then $Z(u)\cong\ZZ/3$. We have $sus\i=u^2$
for some $s\in\D$. We have $s^2us^{-2}=u^4=u$. Thus
$s^2\in Z(u)$. If $s^2\ne1$ then $s$ has order $6$. We
have $s\in Z(s^2)$, a group of order $3$ so $s$ cannot have
order $6$. Thus $s^2=1$.

For any $b\ge1$ we have $su^bs\i=(sus\i)^b=u^{2b}$. 
We have $u^{-b}su^b=u^{-b}u^{2b}s=u^bs$.
Thus $s,us,u^2s$ are all conjugate (and have order $2$).

Let $S$ be the set of (unordered) pairs $\{y,y\i\}$ with
$y\in A_3^8$. This set has $4$ elements. Now $\D$ acts on
$S$ by conjugation. Let $K$ be the set of elements of $\D$
which act trivially on $S$. Let $g\in K$.
We have $Ad(g^2)=1$ on any $H\in S$ hence $H\sub Z(g^2)$.
Then $g^2\in Z(y)$ for any $y\in A_3^8$ hence $g^2=1$ or
$g^2\in\{y,y\i\}$ for any $y\in A_3^8$. (The last
possibility does not arise since there are four sets $H$).
Thus $g^2=1$.
Thus $\sha(K)$ is a sum of some numbers in $\{1,3,6,6\}$
with $1$ definitely appearing. 
We have $1+3+6=10,1+6=7,1+3=4$. Since
$\sha(K)$ must divide $24$ we must have
$\sha(K)=1+3=4$  or $\sha(K)=1$ that is $K=A_8^3\cup\{1\}$
or $K=\{1\}$.
Assume that $A_8^3\sub K$.
Let $g\in A_8^3$. Since $\sha(Z(g))=8$, the order of $g$ is
a power of $2$. If $Ad(g):H@>>>H$ has a fixed point for
some $H\in S$ then $g$ commutes with an element of $A_3^8$
hence $g=1$ (a contradiction) or $g\in A_3^8$ (but it is
also in $A_8^3$ (a contradiction)).
Thus $Ad(g):H@>>>H$ has no fixed point for any $H\in S$
that is $gyg\i=y\i$ for any $y\in A_3^8$.
From $gug\i=u\i$ we obtain $g\in\{s,su,us\}$.
Thus $A_8^3\sub\{s,su,us\}$ hence $A_8^3=\{s,su,us\}$.
Any two distinct elements in $\{s,su,us\}$ have product in
$\{u,u\i\}$; indeed,

$sus=u^2, uss=u, su^2s=uss=u,u^2ss=u^2$,

$usu^2s=uuss=u^2,u^2sus=suus=u.$
\nl
But $\{u,u\i\}$ is not stable by conjugation
($u$ has $8$
conjugates). This is a contradiction. We see that
$A_8^3\not\sub K$. Hence $K=\{1\}$. Hence the homomorphism
$\D@>>>Perm(S)$ with kernel $K$ (given by conjugation) is
injective. Since $\sha(\D)=\sha(Perm(S))=24$, it is an isomorphism, so that $\D\cong S_4$, as desired.
(We have not used 1.3(b),(d).)

We now assume that $c_E=5!$. The conjugacy classes
in $\D$ are
$$A_{120}^1,A_{12}^{10},A_8^{15},A_6^{20},{}'A_6^{20},A_5^{24},A_4^{30}.$$

Let $u\in A_5^{24}$. Then $u$ has order $5$. We can find
$s\in\D$  such that $ sus\i=u^2$.
For $b\ge1$ we have $su^bs\i=(sus\i)^b=u^{2b}$.
Hence $s^au^bs^{-a}=u^{2^ab}$ for any $a\ge0$.
We have $s^2us^{-2}=u^4=u\i$. Thus $s^2\ne1$. We have
$s^4us^{-4}=u^{16}=u$. Thus $s^4\in Z(u)$. If $s^4\ne1$
then it has order $5$.
Since $s^2\ne1$, $s^2$ would have order $10$. We have
$s^2\in Z(s^4)$, a group of order $5$, so $s^2$ cannot
have order $10$. Thus $s^4=1,s^2\ne1$, $s$ has order $4$.
We show:

($*$) If $g\in\D$ has order $4$ then $g\in A_4^{30}$.
(In particular we have $s\in A_4^{30},s\i\in A_4^{30}$).
\nl
We must have $g\in A_{12}^{10}\cup A_8^{15}\cup A_4^{30}$.
If $g\in A_{12}^{10}$ (resp. $g\in A_8^{15}$) then since
$Z(g)\sub Z(g^2)$ we have $\sha(Z(g^2))=12$ or $120$
(resp. $\sha(Z(g^2))=8$ or $120$); but $g^2\ne1$ so
$\sha(Z(g^2))=12$ (resp. $\sha(Z(g^2))=8$) so that
$Z(g)=Z(g^2)$ and $g,g^2$ are conjugate; but $g^2$ has
order $2$, contradiction.. This proves ($*$).

By $(*)$ we can find $w\in\D$ such that $wsw\i=s\i$. Then
$w^2sw^{-2}=ws\i w\i=(wsw\i)\i=s$. Thus $w^2\in Z(s)$ so
that $w^2=s^k$ for some $k$. If $w^2=s$ then
$w$  commutes with $w^2$ hence with $s$ (absurd,
since $wsw\i=s\i$). Similarly $w^2\ne s\i$. Thus $w^2=1$
or $w^2=s^2$. We have $ws^2w\i=(wsw\i)^2=s^{-2}=s^2$.
Thus $w\in Z(s^2)$.
Now for any $j$ we have $s^j\in Z(s^2)$ hence
$Z(s)\sub Z(s^2)$; also $Z(s^2)$ contains $w$ which is not
in $Z(s)$ hence $Z(s^2)$ contains the subgroup

$\{w^js^k;j=0,1,k=0,1,2,3\}$.
\nl
This subgroup has $8$
elements; thus $\sha(Z(s^2))$ is divisible by $8$ so it is
exactly $8$. We see that $s^2\in A_8^{15}$.

Assume that $w^2=s^2$. Then $w$ has order $4$ hence
$w\in A_4^{30}$ hence $w$ is conjugate to $s$ that is
$s=awa\i$ for some $a\in\D$. Then $s^2=aw^2a\i$ that
is $s^2=as^2a\i$ so that $a\in Z(s^2)$. Thus $s,s\i,w,w\i$
are all conjugate in $Z(s^2)$: $s,s\i$ are conjugate by $w$;
$s,w$ are conjugate by $a$; $s\i,w\i$ are conjugate by $a$.
Thus the conjugacy class of $s$ in $Z(s^2)$ has at least
$4$ elements; hence the centralizer of $s$ in $Z(s^2)$ has
at most $8/4=2$ elements. But this centralizer is
$\{1,s,s^2,s^3\}$ so it contains $4$ elements, a contradiction. We see that $w^2=1$.

Now $Z(s^2)=\{1,s,s^2,s\i,w,ws=s\i w, ws^2=s^2w,ws\i=sw\}$.
This is a Weyl group of type $B_2$ with generators $w,sw$. 'Indeed we have

$(sw)^2=swsw=ws\i sw=1,w(sw)w(sw)=wssw=s^2,$

$(sw)w(sw)w=ssww=s^2$.
\nl
The elements $s^2,w,ws,ws^2,ws\i$ of $Z(s^2)$ have order
$2$.

Let $\D'=\D'_E$ be as in 1.3(d).
The $\D'$-coset containing $A_4^{30}$ must be
$A_4^{30}\cup A_6^{20}\cup A_{12}^{10}$ and the other
coset must be
$A_8^{15}\cup {}'A_6^{20}\cup A_5^{24}\cup A_{120}^1$.
Hence this second coset must be $\D'$ and the first one
must be $\D-\D'$. The commutators of $\D$ must be contained
in $\D'$ since $\D/\D'$ is abelian.
Also $A_4^{30}\sub\D-\D'$ hence the elements  of $A_4^{30}$
are not commutators.

We show that $w,ws$ are not conjugate. Assume $ws=xwx\i$
with  $x\in\D$. Then $s=w\i xwx\i$ is a commutator, absurd.

Note that $w,s^2w$ are conjugate (under $s$); also
$ws,ws\i$ are conjugate (under $s$). We see that

$w,sw^2$ are in one conjugacy class; $ws,ws\i$ are in
another conjugacy class.

Now $Z(w)$ contains the subgroup $\{1,w,s^2,ws^2\}$
so that $\sha(Z(w))$ is divisible by $4$;
$Z(ws)$ contains the subgroup $\{1,ws,sw,s^2\}$ so that
$\sha(Z(ws))$  is divisible by $4$.
It follows that one of $w,ws$ is in $A_8^{15}$
and the other is in $A_{12}^{10}$.
Now $(ws)s(ws)\i=s\i$, so, if necessary, we can assume that
$w\in A_8^{15}$, $ws\in A_{12}^{10}$.
The three elements $s^2,w,ws^2$ are all in $A_8^{15}$ and
commute with each other; the product of any two of them is
the third.

They are exactly the elements of $Z(s^2)\cap A_8^{15}$.

Let $S$ be the set consisting of all subsets of
$\D$ of the form $Z(z)\cap A_8^{15}$ for some
$z\in A_8^{15}$. We have an obvious $\D$-equivariant
surjective map $f:A_8^{15}@>>>S$, $z\m Z(z)\cap A_8^{15}$.
We show that the each fibre of $f$ has $3$ elements. It is
enough to show this for the fibre 
at $H=\{s^2,w,ws^2\}$. We have $f(w)=f(ws^2)=f(s^2)$.
Conversely assume that $z\in A_8^{15}$ satisfies $f(z)=H$
that is $Z(z)\cap A_8^{15}=\{s^2,w,ws^2\}$.
Then $z\in Z(s^2)\cap A_8^{15}=\{s^2,w,ws^2\}$.
This proves our claim. It follows that $\sha(S)=15/3=5$.

Let $K$ be the set of all $g$ such that $gHg\i=H$ for any
$H\in S$. Let $g\in K$. Now $Ad(g^6):H\to H$ is $1$ for
all $H$ so that $\cup_HH\sub Z(g^6)$. The union over $H$
is disjoint. (If $h\in H$ then $H=Z(h)\cap A_8^{15}$).
The union has $15$ elements so that $\sha(Z(g^6))\ge15$. It
follows that $g^6=1$. Assume that $g$ has order $6$.
If $Ad(g):H@>>>H$ has a fixed point then g commutes with
some element of $A_8^{15}$  so that it has order dividing
$4$; this is not the case. Thus $Ad(g):H@>>>H$ has no
fixed point so that $Ad(g^3):H@>>>H$ is $1$ for any $H$.
Thus $Z(g^3)$ contains $\cup_HH$ so it contains at least
$15$ elements. But then $g^3=1$, a contradiction. We see
that the order of $g$ is $3,2$ or $1$.
Now $K$ is a union of conjugacy classes of $\D$.
Hence $\sha(K)$ is a sum of numbers in $20,15,10,1$
(with $1$ definitely appearing.) 
Also $\sha(K)$  divides $120$. We have

$20+15+10+1=46,20+10+1=31,20+15+1=36,20+1=21,$

$15+10+1=26,15+1=16,10+1=11$.
\nl
These numbers don't divide $120$. Hence $\sha(K)=1$ so that
$K=\{1\}$.
We see that the homomorphism $\D@>>>Perm(S);g\m(H\m gHg\i)$
is injective. Since these are groups with $120$ elements
this is an isomorphism. We see that $\D\cong S_5$, as
desired. (We did not use 1.3(b).)
This completes the proof of Proposition 1.3.

\subhead 1.5\endsubhead
Fortunately, in 1.3 we cannot have $c_E=6!$. If we did it
would not be clear how to carry out a proof as in 1.4
in that case.

\enddocument